\documentclass[aip,preprint]{revtex4-1}
\usepackage{graphicx}
\usepackage{amssymb}
\usepackage{amsmath}
\usepackage{xurl}
\usepackage{hyperref}
\newtheorem{theorem}{Theorem}
\begin{document}


\title{Denoising-Induced Recovery of Effective Dynamics} 



\author{Narcicegi Kiran}
\email[]{narcicegi.kiran@uni-hamburg.de}

\affiliation{Fachbereich Mathematik, Universität Hamburg, Bundesstraße 55, 20146 Hamburg, Germany}
\affiliation{Instituto de Ci\^encias Matem\'aticas e Computa\c{c}\~ao, Universidade de S\~ao Paulo, S\~ao Carlos, 13566-590, Brazil}

\date{\today}

\begin{abstract}
Recovering dynamical equations from data is a central problem in nonlinear science. In practice,  
measurements are corrupted by noise and often contain fast oscillations that are not part of the slow dynamics of interest. For this reason, data are typically pre-processed before model recovery or classsification. Here, we reveal the effects of filtering, a specific denoising technique, on recovery. 
Our analysis show that equation discovery performed on filtered data recovers, up to a fixed precision, the normal-form (effective) equations rather than the original system.

\end{abstract}


\maketitle 
\begin{quotation}
In theory, uncovering the underlying real-world systems from data can be relatively straightforward when measurements are noise-free. In practice, however, experimental data are rarely clean; consequently, denoising has become a standard procedure. Statistical findings indicate that denoising improves both reconstruction performance and classification accuracy. Nevertheless, this process is often bypassed in equation discovery theory. Here, we reveal the effects of filtering, a specific denoising technique, on recovery. Our analysis show that equation discovery performed on filtered data recovers, up to a fixed precision, the normal-form (effective) equations rather than the original system. 
\end{quotation}
\section{Introduction}

Recovering dynamical equations from data is a central problem in nonlinear science. 
When successful, it allows us to predict future behavior, detect qualitative changes such as bifurcations, and design control strategies based on interpretable models.

In recent years, there has been significant effort devoted to equation discovery. 
Sparse regression methods, dictionary-based approaches, and other machine-learning techniques have been developed to identify governing equations directly from time series data \cite{pereira2025robust, eroglu2020revealing, delabays2025hypergraph, kiran2026unique, li2025novel, tan2023backpropagation, schaeffer2018extracting, donmez2026reconstructing, yang2026reconstruction}. Under ideal conditions, these methods can recover accurate and parsimonious models.

In practice, however, experimental data are rarely clean. 
Measurements are corrupted by noise and often contain fast oscillations that are not part of the slow dynamics of interest. 
For this reason, data are typically pre-processed before model recovery. 
In the field of neuroscience, for instance, it has been reported that applying low-pass filters to signals improves classification accuracy \cite{cardoso2012neuroimaging, kessler2025eeg}.
In experiments involving chemical oscillators, low-pass filters are used to remove high-frequency fluctuations \cite{nijholt2022emergent}. A similar filtering process can also enhance numerical stability in applications within the geosciences \cite{angel2024bathymetry}.

Despite its widespread use, the dynamical impact of filtering remains poorly understood. 
Does filtering simply remove noise, or does it modify the equations that are recovered? 
If the recovered equations differ from the original ones, is there a systematic relationship between them?

The goal of this paper is to clarify these questions. 
We show that filtering can be interpreted as a coordinate transformation related to a normal-form transformation. 
Filtering suppresses rapidly oscillating terms and effectively moves the system toward its averaged dynamics. 
As a result, equation discovery performed on filtered data recovers, up to a fixed precision, the normal-form (averaged) equations rather than the original full system.

We illustrate this mechanism with an example and discuss practical implications for data-driven modeling.

\section{Problem Setup }
We are given noisy time-series data
\begin{equation}
    \big(x^{\text{obs}}(t_k),\,\dot{x}^{\text{obs}}(t_k) \big)_{k=1}^P,
\end{equation}
with \(x^{\text{obs}}(t_k), \dot{x}^{\text{obs}}(t_k) \in \mathbb{R}^N\), generated by observations of a finite-dimensional, deterministic, time-\(T\)-periodic system 
\begin{equation}
        \dot{x} =  F(t,x; \varepsilon)
\end{equation}
where \(\varepsilon\) is a small parameter. We assume that the observed data are generated by this system but are corrupted by measurement noise. Our setup is subject to some technical assumptions listed in the Formal Statements section.

Before performing the equation recovery, the data are passed through a filter: 
\begin{equation}
    y(t) = \mathcal{F}[x^{\text{obs}}](t),
\end{equation}
where $\mathcal{F}$ denotes a smoothing or low-pass filtering operator with \(y(t)\in \mathbb{R}^N\). Its explicit formula is presented in the following section. Our objective is to understand how the transformation $x \mapsto y$ affects the equations recovered from the data.

\subsection{Regression-based Reconstruction}
We briefly introduce the reconstruction procedure.
First, we filter the noisy data using a centered moving average filter of length \(T\), defined by
\begin{equation}
    y(t) = \frac{1}{T}\int_{t-T/2}^{t+T/2}x^{\text{obs}}(s)ds.
\end{equation}

We arrange the filtered data into the following matrices:
\begin{equation}
    Y := 
\begin{pmatrix} y(t_1)^T \\ \vdots \\ y(t_P)^T \end{pmatrix} \in \mathbb{R}^{P \times N}, \:
\dot{Y} := 
\begin{pmatrix} \dot{y}(t_1)^T \\ \vdots \\ \dot{y}^f(t_P)^T \end{pmatrix} \in \mathbb{R}^{P \times N}.
\end{equation}
Let
\begin{equation}
    \Phi:=(\phi_1, \ldots, \phi_M)
\end{equation}
be a given $M$-tuple of sufficiently smooth, real-valued functions. 
We refer to \(\Phi\) as the library, as is standard in the literature. We introduce the library map \(        \Psi: \mathbb{R}^{P\times N} \to \mathbb{R}^{P\times M}\)
defined by
\begin{equation}
  \Psi(Y)=  \begin{pmatrix} \phi_1(t_1,y(t_1)) & \cdots &\phi_M(t_1,y(t_1)) \\ \vdots &\ddots & \vdots \\ \phi_1(t_P,y(t_P)) & \cdots &\phi_M(t_P,y(t_P)) \end{pmatrix}.
\end{equation}

The reconstruction problem is given by
\begin{equation}
    \hat{c}_i =\arg\min_{\tilde{c}_i} \|\Psi(Y) \tilde{c}_i - \dot{Y}_i\|_{*},
\end{equation}
with or without constraints, where \(\|\cdot  \|_{*}\) is any norm on \(\mathbb{R}^P\). We write \(\|\cdot\|_1, \|\cdot\|_2,\) and \(\|\cdot\|_\infty\) for the \(\ell_1,\ell_2\), and \( \ell_\infty\) norms.

\subsection{Averaging Theory} 
We invoke the following result from averaging theory, which will be useful in what follows. When the original system is in its standard form, 
\begin{equation}
        \dot{x}=\varepsilon F^1 (t,x)+ \varepsilon^2F^{[2]}(t, x;\varepsilon), \quad x(0)=a,
    \end{equation}
averaging theory guarantees \cite{sanders2007averaging} the existence of an effective averaged system 
\begin{equation}
    \dot{z} = \varepsilon G(z), \quad z(0) =a,
\end{equation}
where 
\begin{equation}
    G(z) = \frac{1}{T}\int_0^T F^1(s,z)ds,
\end{equation}
such that \(\|x(t)-z(t)\|_1 =\mathcal{O}(\varepsilon)\) for each \(t\in I\) where \(I=[0,L/\varepsilon]\).

\section{Main Results: An Informal Overview}
We first give an informal overview of the main results. Precise statements are given in the Formal Statements section.

\subsection{Filtering as a Coordinate Transformation}
Recall that the filtered data are given by
\begin{equation}
    y(t) = \mathcal{F}[x^{\text{obs}}](t).
\end{equation}
Under suitable assumptions, filtering acts as a near-identity coordinate transformation:
\begin{equation}
y = x + \kappa(t,x),
\end{equation}
where $\kappa$ is small and removes rapidly oscillating components.

\begin{figure*}[ht]
    \centering
\includegraphics[width=.85\textwidth]{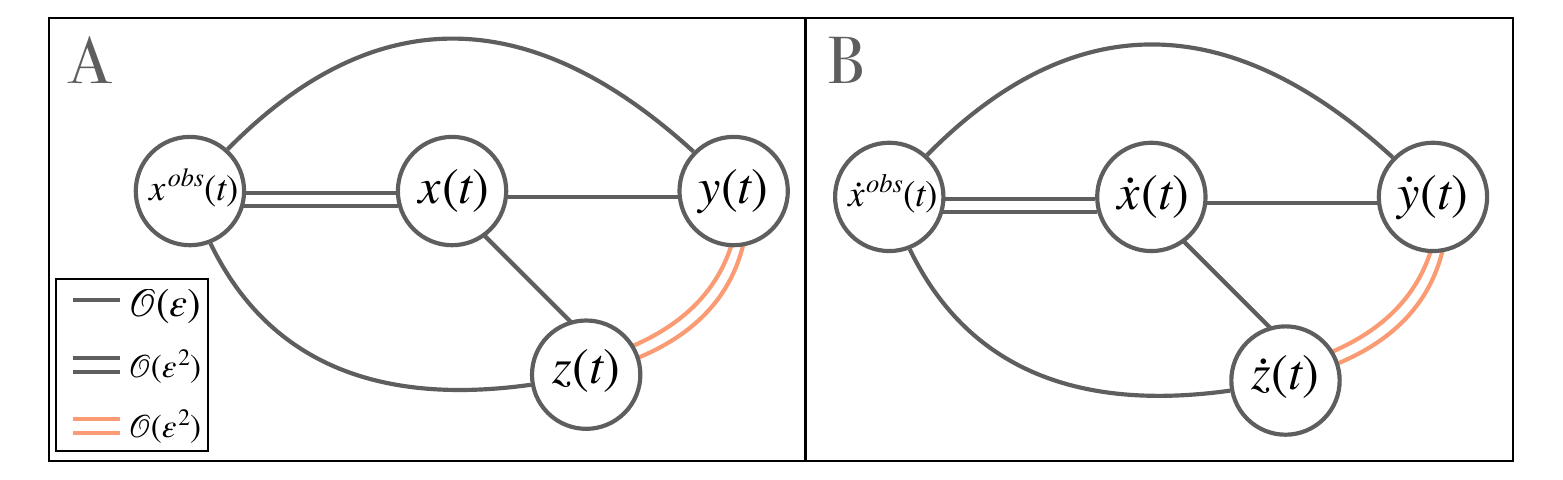}
\caption{Norm estimates leading to the main results. Here, \(x^{\mathrm{obs}}(t), x(t), y(t)\), and \(z(t)\) denote the noisy, original, filtered, and normal-form trajectories, respectively, together with their corresponding derivatives. Estimates shown in grey are already available in the literature, while the estimates shown in orange are proved in this work. All new estimates are in the \(\ell_1\) norm and hold uniformly over the time interval \([T/2, L]\), for some positive constant  \(L\).} 
    \label{fig:estimates}
\end{figure*}

In these new coordinates, the transformed system locally satisfies
\begin{equation}
\dot{y} = G(y) + \text{higher-order terms},
\end{equation}
where $G$ corresponds to the normal-form (averaged) dynamics of the original system. This result follows from norm estimates summarized in Figure \ref{fig:estimates}. Throughout, all \(\mathcal{O}(\cdot)\) statements are as \(\varepsilon\rightarrow 0\), with constants independent of \(\varepsilon\). 
\subsection{Denoising-Induced Reconstruction of Effective Dynamics}

If \(G\) is exactly representable within the given library, that is, there exist coefficients \(c_{ji}\), such that
\begin{equation}
    G_i =\sum_{j=1}^M c_{ji}\phi_j
\end{equation}
for every \(i\in [N]\), then there exists a constructive Basis Pursuit Denoising (BPDN) algorithm that estimates the coefficients of the effective dynamics up to fixed order under mild assumptions.

\section{Illustrative Example}
We illustrate our results with two coupled phase equations that are periodic in time with period \(T\). Coupled phase equations capture the onset of periodic (or almost periodic) behavior, emerging naturally during a phase reduction of coupled dynamics near a Hopf bifurcation. The derivation of the specific form we study, together with the details of our numerical setup, is provided in the Appendix.

We solve the following forward problem 
\begin{equation} \label{system1}
\begin{split}
        \dot{\varphi} &= \varepsilon\big(\delta+\alpha F(t,\varphi,\vartheta)\big),\\
    \dot{\vartheta} &= \varepsilon\big(-\delta+\alpha F(t,\vartheta,\varphi)\big),
\end{split}
\end{equation}
where \begin{equation}
    F(t, \varphi, \vartheta) = \sin(a(\varphi-\vartheta))+cos((b+c)\omega t+b\varphi+c\vartheta)),
\end{equation}
and obtain the noisy observations \(\big( X^{\text{obs}}, \dot{X}^{\text{obs}}\big)\) at a \(\ell_1\)-noise level of \(\mathcal{O}(\varepsilon^2)\). 
We generate filtered observations \(\big( Y, \dot{Y}\big)\) from the noisy data using a centered moving average filter of length \(T\). Assuming the system in Eq. (\ref{system1}) is hidden, we solve BPDN problem using only the filtered observations.

Because fast oscillations are indistinguishable from noise in the original system (especially under \(\mathcal{O}(1/\varepsilon)\) frequency scaling), reconstruction from raw data without filtering is ill-posed.

\begin{figure*}[ht]
    \centering
\includegraphics[width=1\textwidth]{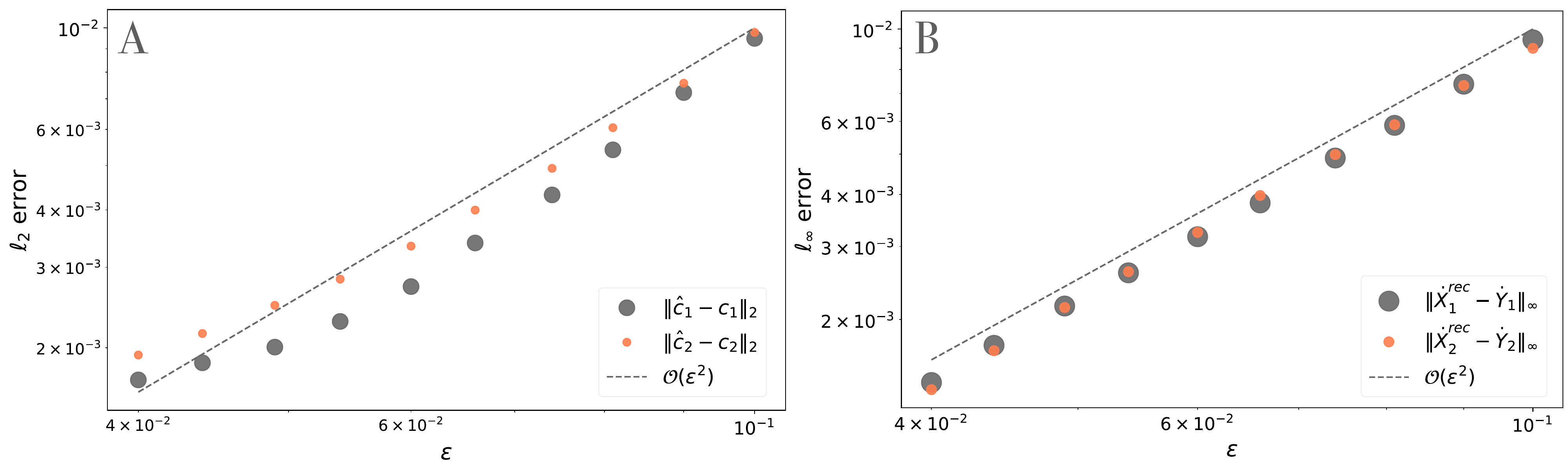}
    \caption{Two coupled phase oscillators. The subscripts 1 and 2 correspond to oscillators \(\varphi\) and \(\vartheta\), respectively. (A) Log–log plot of the \(\ell_2\)-norm diﬀerence
between the estimated coeﬃcient vector \(\hat{c}_i\) and the ground truth coefficient vector of the averaged system \(c_i\) for Basis Pursuit Denoising with tolerance \(\varepsilon^2\), compared against reference rates. (B) Log–
log plot of the sup-norm between the reconstructed vector field and the filtered
derivatives, compared against reference rates. }
    \label{fig:coefficients}
\end{figure*}

In Figure \ref{fig:coefficients}, we show that sparse recovery on filtered data recovers coefficients of the averaged system
\begin{equation}
\begin{split}
        \dot{\bar{\varphi}} &=\varepsilon\big(\delta+\alpha G(\varphi,\vartheta)\big),\\
    \dot{\bar{\vartheta}} &= \varepsilon\big(-\delta+\alpha G(\vartheta,\varphi)\big),
\end{split}
\end{equation}
where \begin{equation}
    G(\varphi, \vartheta) = \sin(a(\varphi-\vartheta))
\end{equation}
with $\ell_2$-error of \( \mathcal{O}(\varepsilon^2)\), which is the order of the observational noise. We note that the library includes all the trigonometric functions appearing above (see Appendix D).

\section{Formal Statements}

Before proceeding to the formal statements, we impose some technical assumptions that are sufficient for our purposes. For a more elaborate version of these results, the reader is referred to the author's thesis \cite{kiran2026recoverability}.

\subsection{Technical Assumptions\label{sec:technical}}

 \textbf{[A1]} Let \(F:\mathbb{R}\times \mathbb{R}^N \times  \mathbb{R}^+\rightarrow \mathbb{R}^N\) be a \(C^r\) vector field with \(r \geq 2\). Consider the dynamical system in standard form
    \begin{equation} \label{standardform}
        \frac{dx}{dt}=F(t,x; \varepsilon), \quad x(t_0) =a,
    \end{equation}
    \(x,a\in U \subseteq \mathbb{R}^N, 0\leq \varepsilon < \varepsilon_0\), with 
    \begin{equation}
        F(t,x;\varepsilon) =\varepsilon F^1 (t,x)+ \varepsilon^2F^{[2]}(t, x;\varepsilon), 
    \end{equation}
    where \(\varepsilon \ll 1\), \(U\) is a bounded, connected, open set, \(F^1\) is Lipschitz continuous, and \(F^{[2]}\) is continuous, both periodic in \(t\) with period \(T\). Here, \({[k]}\) denotes the collection of all higher-order terms of order \(k\) and above. 
    
   \textbf{[A2]} The time-series data are corrupted by observational noise, 
\begin{equation}
\begin{split}
            x^{\text{obs}}(t) = x(t) + \eta^x(t), \\   \dot{x}^{\text{obs}}(t) = \dot{x}(t) + \eta^v(t), 
\end{split}
\end{equation}
    where \(\eta^x(t), \eta^v(t)\in \mathbb{R}^N\) are Riemann integrable satisfying \(\|\eta^x(t)\|_1, \|\eta^v(t)\|_1= \mathcal{O}(\varepsilon^2)\). 
    
 \textbf{[A3]} The data are densely sampled in time, i.e., 
\begin{equation}
    \max_{k}  \Delta t_k \ll \varepsilon,\quad \Delta t_k =t_{k+1}-t_k.
\end{equation}

\textbf{[A4]} We model denoising of noisy data, i.e., \(
    y(t) = \mathcal{F}[x^{\text{obs}}](t)\),  using a centered moving average filter of length \(\xi > \max_k \Delta t_k\), defined for \([t-\xi/2, t+\xi/2]\subset I\) by
\begin{equation}
    y_{\xi}(t) = \frac{1}{\xi}\int_{t-\xi/2}^{t+\xi/2}x^{\text{obs}}(s)ds,
\end{equation}
where \(I =[0,L/\varepsilon]\) for some constant \(L>0\). For notational convenience, we use \(y(t)\) and drop the subscript \(\xi\) when \(\xi =T\).




\subsection{Filtering as a Coordinate Transformation}
The following result enables the interpretation of filtering as a near-identity coordinate transformation.
\begin{theorem}
    Let \(x^{\text{obs}}(t)\) be the solution of the original system in Eq. (\ref{standardform}) corrupted by observational noise satisfying the assumptions above, and let \(y(t)\) denote the filter introduced above. Then, there exists an effective averaged system with trajectory \(z(t)\) such that:
    \begin{enumerate}
        \item The filtered trajectory \(y(t)\) is \(\mathcal{O}(\varepsilon^2)\) close to \(z(t)\), while remaining \(\mathcal{O}(\varepsilon)\) close to the original trajectory subject to observational noise \(x^{\text{obs}}(t)\)  where \(\frac{T}{2}\leq t\leq L\) for some constant \(L>0\).
        \item Similarly, the derivative \(\dot{y}(t)\) is \(\mathcal{O}(\varepsilon^2)\) close to \(\dot{z}(t)\), while remaining \(\mathcal{O}(\varepsilon)\) close to the original derivative subject to observational noise \(\dot{x}^{\text{obs}}(t)\) where \(\frac{T}{2}\leq t\leq L\) for some constant \(L>0\).
    \end{enumerate} 
\end{theorem}

\subsection{Denoising-Induced Reconstruction of Effective
Dynamics} 
The following result formalizes that there exists a constructive Basis Pursuit Denoising (BPDN) algorithm such that the estimated coefficients and the true coefficients of the averaged system remain within the observational noise level. Recall the notions of s-sparse and restricted isometry constant defined in \cite{candes2008restricted}.
\begin{theorem}
    Let the filtered samples \begin{equation}
    \big(y(t_k),\,\dot{y}(t_k) \big)_{k=1}^P,
\end{equation} be given with the technical assumptions above. Assume further that the coefficients of the averaged system \(c_i \in \mathbb{R}^M \) are \(s\)-sparse, and that the restricted isometry constant of library fed by the  averaged trajectories satisfies \(\delta_{2s}(\Psi(Z))< \sqrt{2}-1\). 

Consider the solution \(\hat{c}_i\) of the following minimization problem: 
\begin{equation}
    \hat{c}_i = \arg \min_{\tilde{c}_i \in \mathbb{R}^M}\|\tilde{c}_i\|_1 \text{   subject to   } \|\Psi(Y)\tilde{c}_i -\dot{Y}_i\|_2 \leq k\varepsilon^2
\end{equation}
for some constant \(k>0\). Then there exists a constant \(b>0\) such that \begin{equation}
    \|\tilde{c}_i -c_i\|_2 \leq b \varepsilon^2.
\end{equation}
\end{theorem}
\section{Conclusion}
Here, we present the analytic treatment of the filtering effects on the equation recovery procedure. Recently, significant progress has been made in the reconstruction of normal forms and hyperstructures under conditions of uncertainty and noise \cite{donmez2026reconstructing, nijholt2022emergent, delabays2025hypergraph, yang2026reconstruction}. Filtering not only facilitates the reconstruction of a specific normal form but also offers a means to intervene in a manner sensitive to the noise level. As discussed in the introduction, denoising process yields statistically significant improvements in experimental settings. With this result, we are now able to demonstrate this phenomenon analytically.

\begin{acknowledgments}
We are indebted to Tiago Pereira. We thank Babette de Wolff and Berke Ozgur Arslan for enlightening discussions. We acknowledge support by the Serrapilheira Institute Grant No. 270 Serra-1709-16124. 
\end{acknowledgments}
\section*{Conflict of Interest}
The author has no conflicts to disclose.

\section*{Data Availability}
The data that support the findings of this study are openly available at \url{
https://doi.org/10.5281/zenodo.24222262}.

\appendix
\section{Justification For the Diffusive Coupling in the Illustrative Example} 

Let \(F:\mathbb{T}^2\rightarrow \mathbb{R}\) be a \(C^\infty\) function, and \(2\pi\) periodic function in both variables. Consider two coupled phase variables \(\phi, \theta \in \mathbb{T}\), governed by the equations:
\begin{equation}
    \begin{split}
        \dot{\phi}= \omega_\phi + \varepsilon \alpha F(\phi,\theta),\\
        \dot{\theta}= \omega_\theta + \varepsilon \alpha F(\theta, \phi),
    \end{split}
\end{equation}
\(\varepsilon\ll 1\), and \(\alpha, \omega_\phi, \omega_\theta\) are some constants. Assume that the oscillator frequencies are close to each other and defining reference frequency \(\omega\), we have
\begin{equation}
    \omega_\phi =\omega +\varepsilon\delta_1, \text{ and }\omega_\theta =\omega- \varepsilon\delta_2,
\end{equation}
with \(\delta=\delta_1+\delta_2=\mathcal{O}(1)\). By introducing slow phase variables
\begin{equation}
    \varphi= \phi- \omega t, \text{ and } \vartheta =\theta - \omega t ,
\end{equation}
we obtain
\begin{equation}
\begin{split}
        \dot{\varphi} =\varepsilon \big(\delta_1 + \alpha F(\omega t+ \varphi, \omega t + \vartheta)\big),\\
                \dot{\vartheta} =\varepsilon \big(-\delta_2 + \alpha F(\omega t+ \vartheta, \omega t + \varphi)\big).
\end{split}
\end{equation}
Since \(\varphi, \vartheta \) evolve slowly and \(t\) appears only through the fast phase \(\omega t\), we can average the system over one period of the fast oscillation. Let the averaged interaction be as follows:
\begin{equation}
    G(\varphi,\vartheta) =\frac{1}{2\pi }\int_0^{2\pi} F(s+\varphi, s+\vartheta)ds.
\end{equation}
Then, by transformation \(u = s+ \varphi\), we have 
\begin{equation}
    \begin{split}
        G(\varphi, \vartheta ) = & \frac{1}{2 \pi}\int_{\varphi}^{2\pi + \varphi } F(u, u-\varphi+\vartheta)du\\
        &\frac{1}{2\pi}\int_0^{2\pi} F(u, u +\vartheta-\varphi)du\\
        & =H(\vartheta-\varphi).
    \end{split}
\end{equation}
Therefore, we see that the coupling is diffusive.

\section{Proof of Theorem 1}
We sketch the proof for the filter's window size \(\xi =T\). A similar version of the results can be found in the author's thesis \cite{kiran2026recoverability}.
        
\subsection{Proof of Theorem 1, item (2)}

It is known in the literature\cite{sanders2007averaging} that there exists an invertible near-identity transformation (NIT)
\begin{equation} \label{NIT}
    x(t)=w(t)+\varepsilon u^{[1]}(t,w,\varepsilon),
\end{equation} where \(u^{[1]}\) is periodic in \(t\) with period \(T\) which carries the original system into the full averaged equation 
\begin{equation}
    \dot{w}=\varepsilon G(w) +\varepsilon^2F_*^{[2]}(t,w;\varepsilon), 
\end{equation}
for some \(F_*^{[2]}\), induced by the transformation. We recall that \(F^1\) is \(T\)-periodic in time. So,
\begin{equation}
\begin{split}
        G(x) =& \frac{\varepsilon }{T} \int_{0}^{T} F^1(s,x)ds = \frac{\varepsilon }{T} \int_{l}^{l+T} F^1(s,x)ds,\\
        =&\frac{\varepsilon }{T} \int_{t-T/2}^{t+T/2} F^1(s,x)ds,
\end{split}
\end{equation}
with \(l=t-T/2\). On the other hand, if we filter the system, we have
\begin{equation}
    \begin{split}
      \frac{1}{T} \int_{t-T/2}^{t+T/2}  \dot{x}(s)ds &= \frac{\varepsilon }{T} \int_{t-T/2}^{t+T/2} F^1(s,x)ds +\frac{\varepsilon^2}{T} \int_{t-T/2}^{t+T/2} F_*^{[2]}(s,x;\varepsilon)ds,\\
        &=\varepsilon G(x) +\frac{\varepsilon^2}{T} \int_{t-T/2}^{t+T/2} F_*^{[2]}(s,x;\varepsilon)ds.
    \end{split}
\end{equation}
Another observation is as follows:
\begin{equation}
\begin{split}
        \|x(t)-y(t) \|_1 &\leq \Big\| \frac{1}{T}\int_{t-T/2}^{t+T/2}[x(s,\varepsilon)-x(t,\varepsilon)]ds\Big\|_1 +\frac{1}{T}\int_{t-T/2}^{t+T/2}\|\eta^x(s)\|_1ds,\\
        &\leq  \frac{1}{T}\int_{t-T/2}^{t+T/2}\Big( \int_s^t\|\dot{x}(\tau,\varepsilon)\|_1d\tau\Big) ds +\frac{1}{T}\int_{t-T/2}^{t+T/2}\|\eta^x(s)\|_1ds.
\end{split}
\end{equation}
Since \(\|\dot{x}(\tau,\varepsilon)\|_1 = \mathcal{O}(\varepsilon)\) and \(\|\eta^x(s)\|_1 =\mathcal{O}(\varepsilon^2) \), we have
\begin{equation}
     \|x(t)-y(t) \|_1 = \mathcal{O}(\varepsilon).
\end{equation}
Therefore, 
\begin{equation}
    G(x) = G(y)+\|x-y\|_1 D G(y)\hat{u} +\mathcal{O}(\varepsilon^2),
\end{equation}
with \(\hat{u}=\frac{(x-y)}{\|x-y\|_1}\). Thus,
\begin{equation}
 \dot{y} =\varepsilon G(y) + \mathcal{O}(\varepsilon^2),
\end{equation}
where \( \frac{1}{T} \int_{t-T/2}^{t+T/2}  \dot{x}(s)ds =\dot{y}(t).\) Since \(\dot{z}(t)\) is the truncated version of \(\dot{w}(t)\), \(\|\dot{y}(t)-\dot{z}(t)\|_1 =\mathcal{O}(\varepsilon^2) \), by definition. \(\square\)
\subsection{Proof of Theorem 1, item (1)}

We recall the fundamental theorem of calculus for filtered and truncated averaged systems: 
\begin{equation}
         \begin{split}
              {y}(t)- y(0)&= \int_0^t \dot{y}(s)ds,\\
              z(t)-z(0)&= \int_0^t \dot{z}(s)ds.
         \end{split}
\end{equation}
Therefore, we have
\begin{equation}
    \|z(t)-y(t)\|_1 \leq \|y(0)-z(0) \|_1 + \int_0^t \|\dot{z}(s)-\dot{y}(s)\|_1ds.
\end{equation}
The initial condition term cancels out because \(z(0)=x(0) =a\) and \(y(0) =\frac{1}{T}\int_{t-T/2}^{t+T/2}ad \tau =a\).

 Therefore, there exists a constant \(M>0\) such that 
\begin{equation}
\begin{split}
 \|{z}(t)-y(t)\|_1
        &\leq \int_0^t  \varepsilon\|G(z)-G(y)\|_1ds+ \varepsilon^2Mt, \\
       & \leq \varepsilon \lambda_G \int_0^t \|y(s)-z(s)\|_1ds + \varepsilon^2Mt,
       \end{split}
\end{equation}

where \(\lambda_{G}\) is the Lipschitz constant of \(G\). By Gronwall inequality, we have
\begin{equation}
    \|z(t)-y(t)\|_1 \leq (\epsilon^2Mt) \exp(\varepsilon\lambda_{G}t).
\end{equation}
Thus, for \(t= \mathcal{O}(1)\), we have \(\|y(t) -z(t)\|_1=\mathcal{O}(\varepsilon^2)\).

Lastly, since we defined the time translation \(t=l-T/2\) to apply the filter, the estimates are well-defined after transient time, i.e., for \(t\in[T/2, L]\).
\(\square\)

\section{Proof of Theorem 2}
For notational clarity, we omit the subscript \(i\) in the \(c_i\) and \( \dot{Z}_i, \dot{Y}_i\) terms.

Consider the constraint equation. By the triangle inequality, we have
\begin{equation}
    \|\Psi(Y)\tilde{c}-\dot{Y}\|_2 \leq \|\Psi(Y)\tilde{c}-\Psi(Z)\tilde{c}\|_2+\|\Psi(Z)\tilde{c}-\dot{Y}\|_2.
\end{equation}
By rearranging,
\begin{equation} \label{optimize}
    \|\Psi(Z)\tilde{c}-\dot{Y}\|_2 \leq \|\Psi(Y)\tilde{c}-\Psi(Z)\tilde{c}\|_2+ \|\Psi(Y)\tilde{c}-\dot{Y}\|_2.
\end{equation}
Recall that, by Theorem 1,
\begin{equation}
    \begin{split}
        \|Y-Z\|_2&\leq \sqrt{P}\max_{1\leq  j\leq P} \|y(t_j)-z(t_j)\|_1 =\mathcal{O}(\sqrt{P}\varepsilon^2),\\
        \|\dot{Y}-\dot{Z}\|_2 & \leq  \sqrt{P}\max_{1\leq  j\leq P} \|\dot{y}(t_j)-\dot{z}(t_j)\|_1 =\mathcal{O}(\sqrt{P}\varepsilon^2).
    \end{split}
\end{equation}
Thus, for a fixed number of samples \(P\),
\begin{equation}
    \begin{split}
        Y=&Z+\xi \text{ with } \|\xi\|_2 =\mathcal{O}(\varepsilon^2),\\
        \dot{Y}=&\dot{Z}+r \text{ with } \|r\|_2 =\mathcal{O}(\varepsilon^2).
    \end{split}
\end{equation}
We now consider the first term on the right-hand side in equation \ref{optimize}. Since the library map \(\Psi\) elements are smooth and the admissible set is compact by construction, their first derivatives are uniformly bounded. Therefore, 
\begin{equation}
    \|\Psi(Y)-\Psi(Z)\|_{op} \leq b_1 \|Y-Z\|_2 =\mathcal{O}(\varepsilon^2),
\end{equation}
for some constant \(b_1>0\). Consequently, 
\begin{equation}
    \|\Psi(Y)\tilde{c}-\Psi(Z)\tilde{c}\|_2 \leq b_1 \varepsilon^2 \|\tilde{c}\|_2.
\end{equation}
Combining this with feasibility constraint 
\begin{equation}
    \|\Psi(Y)\tilde{c}-\dot{Y}\|_2 \leq k\varepsilon^2,
\end{equation}
we obtain
\begin{equation}
    \|\Psi(Z)\tilde{c}- \dot{Y}\|_2 \leq k\varepsilon^2 + b_1 \varepsilon^2\|\tilde{c}\|_2.
\end{equation}
Thus, optimization problem reduces to the following problem:
\begin{equation}
    \hat{c} =\arg \min_{\tilde{c}\in \mathbb{R}^M}\|\tilde{c}\|_1 \text{ subject to } \|\Psi(Z)\tilde{c}-\dot{Y}\|\leq b_2 \varepsilon^2,
\end{equation}
for some constant \(b_2>0\). This is exactly the noisy measurement setting in \cite{candes2008restricted}, and the RIP condition guarantees that the solution satisfies 
\begin{equation}
    \|\hat{c}-c\|_2 \leq b\varepsilon^2,
\end{equation}
for some constant \(b>0\).

\section{Numerical Details}

The Python code is publicly available on Zenodo \cite{kiran2026filter_code}. The details generating the results for the two coupled phase oscillators in Figure \ref{fig:coefficients} are as follows:

\begin{itemize}
    \item All numerical illustrations were performed using representative parameters \(\omega=3, \delta=1.2,\alpha=0.9, a=3, b=10,c=15\) and matched initial conditions for the original and averaged systems, ensuring consistent comparison of trajectories and coeﬃcient recovery.
    \item For each \(\varepsilon\) in the prescribed range, we generate trajectories \(X=(\varphi, \vartheta)\) and their derivatives by solving the system and sampling the solution on a uniform time grid, and add Gaussian noise of amplitude \(\mathcal{O}(\varepsilon^2)\) to both the state variables and their time derivatives.
    \item The noisy trajectories and derivatives are filtered using a time-based moving average
with window size \(\pi/\omega\).
\item Using the filtered trajectories sampled on a uniform time grid, we construct the library
and solve a Basis Pursuit Denoising problem on a masked time interval to mitigate filtering-induced boundary eﬀects, using the filtered derivatives as observations where the library is given as
\begin{equation}
\begin{split}
        \Phi =( 1, & \sin(3(\varphi-\vartheta)), \cos(3(\varphi-\vartheta)), \sin(2(\varphi+\vartheta)), \\&\cos(2(\varphi+\vartheta)), \sin(10\varphi+15\vartheta+25t),\cos(15\varphi+10\vartheta+25t) ),
\end{split}
\end{equation}
where low-frequency terms capture the averaged dynamics, high-frequency terms capture fast oscillations. 
\item As we claimed, we also generate the averaged system’s trajectories and derivatives,
i.e., reduced reference system. 
\item To illustrate our statement, we compute the \(\ell_2\)-error between the recovered coeﬃcient
vector and the ground-truth coeﬃcients corresponding to the averaged system, and
compare it against the reference line \(\mathcal{O}(\varepsilon^2)\).
\item Furthermore, we plot the supremum norm of the reconstruction error in the time
derivatives as a function of \(\varepsilon\), and compare it against the same reference rates. 
\end{itemize}
We remark that the RIP condition is not satisfied, but the columns sampled along the trajectory are eﬀectively independent due to ergodicity, enabling
accurate sparse recovery with \(\mathcal{O}(\varepsilon^2)\) error.

%
%

%


\bibliography{denoising}

\end{document}